\documentclass[12pt]{amsart}
\usepackage{amsfonts}
\usepackage{amsmath,amssymb,amsthm}

\usepackage{color}
\newcommand{\R}{\mathbb R}

\newcommand{\ds}{\displaystyle}
\newtheorem{thm}{Theorem}[section]
\newtheorem{cor}[thm]{Corollary}
\newtheorem{lem}[thm]{Lemma}
\newtheorem{prop}[thm]{Proposition}
\theoremstyle{definition}
\newtheorem{defn}[thm]{Definition}
\theoremstyle{remark}
\newtheorem{rem}[thm]{Remark}

\newcommand\topstrut{\rule{0mm}{2.9ex}}
\newcommand\bottomstrut{\rule[-1.5ex]{0mm}{1.5ex}}
\newcommand\titlestrut{\topstrut\bottomstrut}
\begin{document}

\title[SPACE-LIKE WEINGARTEN SURFACES IN MINKOWSKI SPACE]
{Space-like Weingarten surfaces in the three-dimensional Minkowski Space
and their natural Partial Differential Equations}%

\author{Georgi Ganchev and Vesselka Mihova}%
\address{Bulgarian Academy of Sciences, Institute of Mathematics and Informatics,
Acad. G. Bonchev Str. bl. 8, 1113 Sofia, Bulgaria}%
\email{ganchev@math.bas.bg}%
\address{Faculty of Mathematics and Informatics, University of Sofia,
J. Bouchier Str. 5, 1164 Sofia, Bulgaria}
\email{mihova@fmi.uni-sofia.bg}


\begin{abstract}

On any space-like Weingarten surface in the three-dimensional
Minkowski space we introduce locally natural principal parameters
and prove that such a surface is determined uniquely up to motion by
a special invariant function, which satisfies a natural non-linear
partial differential equation. This result can be interpreted as a
solution to the Lund-Regge reduction problem for space-like
Weingarten surfaces in Minkowski space. We apply this theory to
linear fractional space-like Weingarten surfaces and obtain the
natural non-linear partial differential equations describing them.
We obtain a characterization of space-like surfaces, whose
curvatures satisfy a linear relation, by means of their natural
partial differential equations. We obtain the ten natural PDE's
describing all linear fractional space-like Weingarten surfaces.
\end{abstract}

\maketitle

\section{Introduction}
It has been known to Weingarten \cite{W1, W2}, Eisenhart \cite{E},
Wu \cite{Wu}) that without changing the principal lines on a
Weingarten surface, one can find geometric coordinates in which the
coefficients of the metric are expressed by the principal curvatures
(or radii of curvature).

The geometric parameters on Weingarten surfaces were used in
\cite{Wu} to find the classes of Weingarten surfaces yielding
``geometric" $\mathfrak{so}(3)$-scattering systems (real or complex)
for the partial differential equations, describing these surfaces.

In \cite{GM1} we have shown that the Weingarten surfaces
(W-surfaces) in Euclidean space admit geometrically determined
principal parameters (\emph{natural principal parameters}), which
have the following property: all invariant functions on W-surfaces
can be expressed in terms of one function $\nu$, which satisfies one
\emph{natural} partial differential equation. The Bonnet type
fundamental theorem states that any solution to the natural partial
differential equation determines a W-surface uniquely up to motion.

Thus the description of any class of W-surfaces (determined by a
given Weingarten relation) is equivalent to the study of the
solution space of their natural PDE. This solves the Lund-Regge
reduction problem \cite{LR} for W-surfaces in Euclidean space.

In this paper we study space-like surfaces in the three dimensional Minkowski
space $\R^3_1$.

A space-like surface $S$ with principal normal curvatures $\nu_1$
and $\nu_2$ is a Weingarten surface (W-surface) \cite{W1, W2} if there exists a
function $\nu$ on $S$ and two functions (Weingarten functions) $f, \, g$ of one
variable, such that
$$\nu_1=f(\nu), \quad \nu_2=g(\nu).$$

A basic property of W-surfaces in Euclidean space is the following theorem of Lie \cite{Lie}:

\emph{The lines of curvature of any W-surface can be found in quadratures.}

This remarkable property is also valid for space-like W-surfaces in Minkowski space.

We use four invariant functions (two principal normal curvatures $\nu_1,\,\nu_2$
and two principal geodesic curvatures $\gamma_1,\, \gamma_2$) and divide space-like
W-surfaces into two classes with respect to these invariants:

(1) the class of \emph{strongly regular} space-like surfaces defined by
$$(\nu_1-\nu_2)\,\gamma_1 \,\gamma_2\neq 0;$$

(2) the class of space-like surfaces defined by
$$\gamma_1= 0, \quad (\nu_1-\nu_2)\, \gamma_2\neq 0.$$

The basic tool to investigate the relation between space-like surfaces and
the partial differential equations describing them, is Theorem 2.1. This theorem
is a reformulation of the fundamental Bonnet theorem for the class of strongly
regular space-like surfaces in terms of the four invariant functions.
Further, we apply this theorem to space-like W-surfaces.

In Section 3 we prove (Proposition 3.4) that any space-like W-surface admits
locally special principal parameters (\emph{natural principal parameters}).

Theorem 3.7 is the basic theorem for space-like W-surfaces of type
(1):

\emph{Any strongly regular space-like W-surface is determined
uniquely up to motion by the functions $f$, $g$ and the function
$\nu$, satisfying the natural PDE $(3.4)$.}

Theorem 3.8 is the basic theorem for space-like W-surfaces of type
(2):

\emph{Any space-like W-surface with $\gamma_1=0$ is determined
uniquely up to motion by the functions $f$, $g$ and the function
$\nu$, satisfying the natural ODE $(3.6)$.}

In natural principal parameters the four basic invariant functions, which determine
space-like W-surfaces uniquely up to motions in $\R^3_1$, are expressed by a single
function, and the system of Gauss-Codazzi equations reduces to a single partial
differential equation (the Gauss equation). Thus, the number of the four invariant
functions, which determine space-like W-surfaces, reduces to one
invariant function, and the number of Gauss-Codazzi equations reduces to one
\emph{natural} PDE. This result gives a solution to the Lund-Regge reduction
problem \cite{LR} for the space-like W-surfaces in $\R^3_1$. The Lund-Regge
reduction problem has been analyzed and discussed from several view points in
the paper of Sym \cite{Sym}.

The class of space-like W-surfaces contains all classical special surfaces:
maximal space-like surfaces (with zero mean curvature $H=0$), space-like
CMC-surfaces (with constant mean curvature $H={\rm const} \neq 0$), space-like
surfaces of constant Gauss curvature $K=\pm 1$ etc. All these surfaces are
examples of space-like W-surfaces with simple Weingarten functions.

A more general class, which contains the above mentioned surfaces, is the class of
linear fractional space-like W-surfaces defined by the following relation between
the normal curvatures $\nu_1$ and $\nu_2$:
$$\nu_1=\frac{A\nu_2+B}{C\nu_2+D}\,, \quad A,B,C,D - {\rm constants}; \quad
BC-AD\neq 0.$$

We prove that on any linear fractional space-like W-surface the
extrinsic Gauss curvature $K'=\nu_1\,\nu_2$, the mean curvature $H$
and the curvature $H'=\ds{\frac{\nu_1-\nu_2}{2}}$  satisfy the
linear relation
$$\delta K' = \alpha H + \beta H' + \gamma, \quad \alpha, \beta, \gamma, \delta -
{\rm constants}; \quad \alpha^2-\beta^2 + 4 \gamma\delta \neq 0,\leqno(1.1)$$
and vice versa (cf \cite{GM2}).

On space-like surfaces in the three-dimensional Minkowski space the
(intrinsic) Gauss curvature $K$ and the extrinsic Gauss curvature
$K'$ are in the relation $K'=-K$.

We note that the class of linear fractional space-like surfaces is invariant under
parallel transformations of space-like surfaces.

In Proposition 4.1 we prove that {\it the natural principal parameters of a given space-like
W-surface $S$ are natural principal parameters for all parallel space-like surfaces $\bar S(a),\,
a={\rm const}\neq 0$ of $S$}.

In Theorem 4.2 we prove the following essential property of the space-like surfaces, which
are parallel to a given space-like surface (cf \cite{GM2}):

{\it The natural PDE of a given space-like W-surface $S$ is the natural PDE of any parallel
space-like surface $\bar S(a),\; a={\rm const}\neq 0$, of $S$.}

Using this theorem, in Section 5 we classify the natural PDE's of the linear fractional
space-like W-surfaces, or equivalently the space-like surfaces satisfying the relation (1.1).

To illustrate our investigations, let us consider the space-like surfaces in Minkowski
space, which satisfy the linear relation
$$\delta K'=\alpha H + \gamma, \qquad \alpha^2+4\gamma\delta \neq 0.\leqno(1.2)$$

Such surfaces are studied in \cite{Mil, Bu, FG, Fo} (see also
\cite{C}).

All surfaces, satisfying the above linear relation, are linear fractional space-like
W-surfaces and each of them has a natural PDE.

Milnor has proved in \cite {Mil} that any space-like surface in $\R^3_1$ satisfying
(1.2) is parallel to a surface, satisfying one of the following conditions:
$H=0$, $K=1$ or $K=-1$.

Here we prove the following result:
\begin{itemize}
\item The space-like surfaces, which are parallel to space-like surfaces with $H=0$,
are described by the natural PDE
$$ \lambda_{xx} + \lambda_{yy} = e^{\lambda}.$$

\item The space-like surfaces, which are parallel to space-like surfaces with
$H = {\rm const}$, are described by the one-parameter family of natural PDE's
$$\lambda_{xx} + \lambda_{yy} = 2|H| \sinh \lambda.$$

Up to similarity, the space-like surfaces, which are parallel to space-like surfaces
with
$H = {\rm const}$, are described by the natural PDE of the surfaces with $H=1/2$\,.
\vskip 2mm

\item The space-like surfaces, which are parallel to space-like surfaces with
$K={\rm const} > 0$ $\,(K'={\rm const}<0)$, are described by the one-parameter family
of natural PDE's
$$ \lambda_{xx} - \lambda_{yy} = -K^2 \sin  \lambda.$$

Up to similarity, the space-like surfaces, which are parallel to space-like surfaces with
$K={\rm const} > 0$, are described by the natural PDE of the surfaces with $K = 1$ $\,(K'=-1)$.
\end{itemize}

We call the surfaces with $H=0$, $H=1/2$ and $K=-K'=1$ \emph{the basic surfaces} in the
class of surfaces, determined by the relation (1.2).

Then we have:

{\it Up to similarity, the surfaces, whose curvatures satisfy the linear relation $(1.2)$,
are described by the natural PDE's of the basic surfaces.}

We apply the above scheme to obtain all natural PDE's describing the space-like surfaces,
whose curvatures $K'$, $H$ and $H'$ satisfy the linear relation (1.1)

We divide the class of all linear fractional space-like W-surfaces
into ten subclasses, which are invariant under parallel
transformations. For any of these  classes we determine the
\emph{basic surfaces} in terms of suitable constants $p$ and $q$.
Applying Theorem 4.2, we find the PDE's which describe all
space-like surfaces in Minkowski space, whose curvatures satisfy the
linear relation (1.1) (Theorem 5.2).

It is essential to note that the natural PDE's of the linear
fractional space-like W-surfaces are expressed in the form $ \delta
\lambda = f(\lambda)\,$, \;where $\delta$ is one of the operators
(cf \cite{GM2}):

$$ \Delta \lambda= \lambda_{uu} + \lambda_{vv},\qquad
\bar {\Delta} \lambda= \lambda_{uu} - \lambda_{vv},$$

$$ \Delta^* \lambda= \lambda_{uu} + \left(\frac{1}
{\lambda}\right)_{vv},\qquad
 \bar{\Delta}^* \lambda= \lambda_{uu} - \left(\frac{1}
{\lambda}\right)_{vv}.$$

\vskip 2mm
In the Euclidean space the intrinsic Gauss curvature $K$ and the extrinsic Gauss curvature $K'$ coincide.

Next we draw a parallel between the natural PDE's of the basic classes of
linear fractional W-surfaces in Euclidean space \cite{GM2} and the natural PDE's of the
corresponding classes of linear fractional space-like surfaces in Minkowski space:
\vskip 6mm

\begin{tabular}{|c|c|c|c|}
\hline
\titlestrut
Nr & Basic surfaces &
Natural PDE  & Natural PDE\\

& \titlestrut $p$, $q$ - {\rm constants} & in $\mathbb{R}^3$  &
in $\mathbb{R}_1^3$\\
\hline
\titlestrut

1 & $H=0$ &$\Delta \lambda = - e^{\lambda}$ & $\Delta \lambda =  e^{\lambda}$\\
[4.0ex]

2 & $H=\frac{1}{2}$ & $\Delta \lambda = -\sinh{\lambda}$ & $\Delta \lambda =  \sinh{\lambda}$\\
[4.0ex]

3 & $H'=1$& $\Delta^*(e^{\nu})=-2\,\nu\,(\nu+2)$ &
$\Delta^*(e^{\nu})=2\,\nu\,(\nu+2)$\\
[4.0ex]

4 & $\begin{array}{c} H=p\,H',\\
p^2>1\end{array}$&$\Delta^*(\nu^{p})=
-2\,\frac{p\,(p+1)}{(p-1)^2}\,\nu$& $\Delta^*(\nu^{p})=
2\,\frac{p\,(p+1)}{(p-1)^2}\,\nu$\\
[4.0ex]

5 & $\begin{array}{c} H=p\,H',\\
p^2<1,\;p\neq 0\end{array}$&$\bar\Delta^*(\nu^{p})=
-2\,\frac{p\,(p+1)}{(p-1)^2}\,\nu$& $\bar\Delta^*(\nu^{p})=
2\,\frac{p\,(p+1)}{(p-1)^2}\,\nu$\\
[4.0ex]

6 & $\begin{array}{c} H=p\,H'+1,\\
p^2>1\end{array}$&$\Delta^*(\lambda^{p})=
-\frac{p\,((p-1)\lambda+2)((p+1)\lambda+2)}
{2\,(p-1)\,\lambda}$&
$\Delta^*(\lambda^{p})=
\frac{p\,((p-1)\lambda+2)((p+1)\lambda+2)}
{2\,(p-1)\,\lambda}$  \\
[4.5ex]

7 & $\begin{array}{c} H=p\,H'+1,\\
p^2<1,\; p\neq 0\end{array}$&$\bar\Delta^*(\lambda^{p})=
-\frac{p\,((p-1)\lambda+2)((p+1)\lambda+2)}
{2\,(p-1)\,\lambda}$&
$\bar\Delta^*(\lambda^{p})=
\frac{p\,((p-1)\lambda+2)((p+1)\lambda+2)}
{2\,(p-1)\,\lambda}$  \\
[4.5ex]

8& $K'=-1$ & $\bar\Delta \lambda = \sin \lambda$  &$\bar\Delta \lambda =- \sin \lambda$\\
[4.5ex]

9 & $ K'=2\,H'$& $\Delta^*(e^{\lambda})=-2
$&
$\Delta^*(e^{\lambda})=2$\\
[4.5ex]

10 & $\begin{array}{c} K'=p\,H'-q,\\
p\neq 0,\,q>0\end{array}$& $\begin{array}{c}
\Delta^*(e^{p\,\mathcal{I}}) =-\frac{p\,q}{2}\,
\frac{\lambda\,\left(p\,\lambda-2\,q\right)}{\lambda^2+q},\\
\mathcal{I}=\frac{1}{\sqrt{q}}\,\arctan\frac{\lambda}{\sqrt{q}}
\end{array}$&
$\begin{array}{c}\Delta^*(e^{p\,\mathcal{I}}) =\frac{p\,q}{2}\,
\frac{\lambda\,\left(p\,\lambda-2\,q\right)}{\lambda^2+q},\\
\mathcal{I}=\frac{1}{\sqrt{q}}\,\arctan\frac{\lambda}{\sqrt{q}}\end{array}$\\
[4.0ex]
\hline
\end{tabular}
\vskip 6mm

{\it Let $S$ be a linear fractional W-surface in $\R^3$ and $S'$ be its corresponding
space-like linear fractional W-surface in $\R^3_1$ generated by the same Weingarten
functions as $S$. Then the natural PDE's of $S$ and $S'$ are related as follows:
$$\R^3:\; \delta \lambda = f(\lambda) \qquad \iff \qquad
\R^3_1: \; \delta \lambda = - f(\lambda),$$
where $\delta$ is one of the operators $\Delta,\, \bar{\Delta}, \,  \Delta^*,$
and \,$\bar \Delta^*$.}
\vskip 4mm

Hu elucidated in \cite{Hu1} the relationship between the PDE's
$$\alpha_{uu}-\alpha_{vv}=\pm\sin \alpha \quad(\sin{\rm -Gordon \; PDE}),$$
$$\alpha_{uu}-\alpha_{vv}=\pm\sinh \alpha \quad(\sinh{\rm -Gordon \; PDE}),$$
$$\alpha_{uu}+\alpha_{vv}=\pm\sin \alpha \quad(\sin{\rm -Laplace \; PDE}), $$
$$\alpha_{uu}+\alpha_{vv}=\pm\sinh \alpha \quad(\sinh{\rm -Laplace \; PDE})$$
and the construction of various kinds of surfaces of constant curvature in $\mathbb{R}^3$
or $\mathbb{R}^3_1$.
\vskip 2mm

The result of Baran and Marvan in \cite{BM1} asserts that "the
simple relation $\rho_1-\rho_2= {\rm const}$ between the principal
radii of curvature determines an integrable class of Weingarten
surfaces in the Euclidean space (the forgotten class)" (see also
\cite{GM1, Rib, vL1, vL2}). The authors associate with the
"forgotten class" the nonlinear partial differential equation
$\displaystyle{z_{uu}+\left(\frac{1}{z}\right)_{vv}+2=0}$ and prove
that it is integrable. These surfaces correspond to the class 9 in
the above table.

In \cite{BM2} they give classes of Weingarten surfaces in Euclidean
space integrable in the sense of soliton theory.

In \cite{Hu2} by using Darboux transformations, from a known solution to
the $\sinh$-Laplace (resp. $\sin$-Laplace) equation  have been obtained explicitly
new solutions of the $\sin$-Laplace (resp. $\sinh$-Laplace) equation.

It is interesting to find the connection between parallel
transformations of pseudo-spherical surfaces ($K = -1$) and their
B\"acklund transformations. A nice base for further investigations
in this direction gives the paper of  S. Buyske \cite{Bu}. This
paper shows that many of the problems in Euclidean space have their
exact analogues in Minkowski three space.

\section{Preliminaries}

In this section we introduce the denotations and formulas in the
theory of surfaces in Minkowski space, which we use further.

Let $\R^3_1$ be the three dimensional Minkowski space with the standard flat
metric $\langle \, , \, \rangle$ of signature $(2,1)$. We assume that the
following orthonormal coordinate system $Oe_1e_2e_3: \; e_1^2=e_2^2=-e_3^2=1, \;
e_i\,e_j=0, \, i\neq j$ is fixed and gives the orientation of the space.

Let $S: \, z=z(u,v), \; (u,v) \in {\mathcal D}$\, be a space-like surface in
${\R}_1^3$ and $\nabla$ be the flat Levi-Civita
connection of the metric $\langle \, , \, \rangle$. The unit normal vector field
to $S$ is denoted by $l$ and $E, F, G; \; L, M, N$ stand for the coefficients of
the first and the second fundamental forms, respectively. Then we have
$$E=z_u^2 >0, \quad F=z_u z_v,\quad G=z_v^2 >0, \quad  EG-F^2>0, \quad l^2=-1;$$
$$L=lz_{uu},\quad M=lz_{uv},\quad N=lz_{vv}.$$

We suppose that the surface has no umbilical points and the principal lines on
$S$ form a parametric net, i.e.
$$F(u,v)=M(u,v)=0, \quad (u,v) \in \mathcal D.$$
Then the principal curvatures $\nu_1, \nu_2$ and the principal geodesic curvatures
(geodesic curvatures of the principal lines) $\gamma_1, \gamma_2$ are given by
$$\nu_1=\frac{L}{E}, \quad \nu_2=\frac{N}{G}; \qquad
\gamma_1=-\frac{E_v}{2E\sqrt G}, \quad \gamma_2= \frac{G_u}{2G\sqrt E}. \leqno(2.1)$$

The mean curvature, the (intrinsic) Gauss curvature and the extrinsic Gauss curvature of $S$
are denoted by $H$, $K$ and $K'$, respectively.

We consider the tangential frame field $\{X, Y\}$ defined by
$$X:=\frac {z_u}{\sqrt E}, \qquad Y:=\frac{z_v}{\sqrt G}$$
and suppose that the moving frame $XYl$ is always positive oriented.

Without loss of generality we consider space-like surfaces satisfying the condition
(cf \cite{GM1})
$$\nu_1(u,v)-\nu_2(u,v)>0, \qquad (u,v) \in \mathcal D. \leqno(2.2)$$

The frame field $XYl$ satisfies the following Frenet type formulas:

$$\begin{tabular}{ll}
$\left|\begin{array}{llccc}
\nabla_{X} \,X & = &  &\gamma_1 \,Y - \nu_1 \, l,  &\\
[2mm]
\nabla_{X} Y & = -\gamma_1 \, X, & & \\
[2mm]
\nabla_{X} \, l & = - \nu_1 \, X; & & &
\end{array}\right.$ &
\quad
$\left|\begin{array}{llccc}
\nabla_{Y} \,X & = & & \gamma_2 \, Y, &\\
[2mm]
\nabla_{Y} Y & = -\gamma_2 \, X & &  & - \nu_2 \, l,\\
[2mm]
\nabla_{Y}\, l & = &  & - \nu_2 \, Y. &
\end{array}\right.$
\end{tabular}\leqno(2.3)$$
\vskip 2mm

The Codazzi equations have the form
$$\gamma_1=\frac{Y(\nu_1)}{\nu_1-\nu_2}=
\frac{(\nu_1)_v}{\sqrt G\,(\nu_1-\nu_2)}\,, \qquad
\gamma_2=\frac{X(\nu_2)}{\nu_1-\nu_2}=
\frac{(\nu_2)_u}{\sqrt E\,(\nu_1-\nu_2)}\,,
\leqno(2.4)$$
and the Gauss equation can be written as follows:
$$K=Y(\gamma_1)-X(\gamma_2)- (\gamma_1^2+\gamma_2^2) = -\nu_1\nu_2= -K',$$
or
$$\frac{(\gamma_2)_u}{\sqrt E}-\frac{(\gamma_1)_v}{\sqrt G}+
\gamma_1^2+\gamma_2^2=\nu_1\nu_2= K'.\leqno(2.5)$$

The Codazzi equations (2.4) imply the following equivalence
$$\gamma_1 \gamma_2 \neq 0 \; \iff \; (\nu_1)_v (\nu_2)_u \neq 0.$$

We consider two types of space-like surfaces parameterized by principal parameters:

\begin{itemize}
\item  \emph{strongly regular} space-like surfaces determined by the condition (cf \cite{GM1})
$$(\nu_1-\nu_2)\gamma_1(u,v)\gamma_2(u,v) \neq 0, \quad (u,v) \in \mathcal D,$$

\item  space-like surfaces satisfying the conditions
$$\gamma_1(u,v)= 0,\quad (\nu_1-\nu_2)\gamma_2(u,v)\neq 0, \quad (u,v) \in \mathcal D.$$
\end{itemize}

Because of (2.4) the formulas
$$\sqrt E=\frac{(\nu_2)_u}{\gamma_2\,(\nu_1-\nu_2)} >0, \quad
\sqrt G=\frac{(\nu_1)_v}{\gamma_1\,(\nu_1-\nu_2)}>0 \, \leqno(2.6)$$
are valid on strongly regular space-like surfaces.

Taking into account (2.6), for strongly regular space-like surfaces formulas (2.3) become

$$\left|\begin{array}{l}
\displaystyle{X_u  =  \frac{\gamma_1 \, (\nu_2)_u}
{\gamma_2(\nu_1-\nu_2)}\, Y- \frac{\nu_1 \, (\nu_2)_u}{\gamma_2
(\nu_1-\nu_2)}\, l,\; Y_u  =  - \frac{\gamma_1 \,
(\nu_2)_u} {\gamma_2(\nu_1-\nu_2)}\, X,\;  l_u  =   -
\frac{\nu_1 \, (\nu_2)_u}{\gamma_2(\nu_1-\nu_2)}\, X;}\\
[4mm]
\displaystyle{X_v =\frac{\gamma_2 \, (\nu_1)_v}
{\gamma_1(\nu_1-\nu_2)}\, Y, \; Y_v  = -\frac{\gamma_2
\, (\nu_1)_v} {\gamma_1(\nu_1-\nu_2)}\, X   -
\frac{\nu_2 \,(\nu_1)_v}{\gamma_1 (\nu_1-\nu_2)}\, l,\;
 l_v  =  -\frac{\nu_2 \, (\nu_1)_v}
{\gamma_1(\nu_1-\nu_2)}\,Y.}
\end{array}\right.\leqno(2.7)$$

Now, finding the compatibility conditions for the system (2.7), we can
reformulate the fundamental Bonnet theorem for strongly regular space-like
surfaces in ${\R}^3_1$ in terms of the invariants of the surface.

\begin{thm}\label{T:2.1}
Let the four functions $\nu_1(u,v), \, \nu_2(u,v), \, \gamma_1(u,v),
\, \gamma_2(u,v)$ be defined in a neighborhood \,$\mathcal D$ of $(u_0, v_0)$,
and satisfy the conditions
$$\begin{array}{ll}
1) & \gamma_1\,(\nu_1)_v >0, \qquad \gamma_2\,(\nu_2)_u >0; \\
[3mm]
2.1) & \ds{\left(\ln \frac{(\nu_1)_v}{\gamma_1}\right)_u=
\frac{(\nu_1)_u}{\nu_1-\nu_2}\,,
\quad \left(\ln \frac{(\nu_2)_u}{\gamma_2}\right)_v=
-\frac{(\nu_2)_v}{\nu_1-\nu_2}\,;}\\
[6mm]
2.2) & \displaystyle{\frac{\nu_1-\nu_2}{2}
\left(\frac{(\gamma_2^2)_u}{(\nu_2)_u}-
\frac{(\gamma_1^2)_v}{(\nu_1)_v}\right)+ \gamma_1^2+\gamma_2^2=
\nu_1\,\nu_2}\,.
\end{array}$$
Let $ z_0\,X_0\,Y_0\,l_0$ be an initial positive oriented orthonormal frame.

Then there exists a unique up to a motion strongly regular space-like surface
$S: \; z=z(u,v), \; (u,v) \in \mathcal D_0 \;
((u_0, v_0) \in \mathcal D_0 \subset \mathcal D)$ with prescribed invariants
$\nu_1, \, \nu_2, \, \gamma_1, \, \gamma_2$ such that
$$z(u_0, v_0)=z_0, \; X(u_0, v_0)=X_0, \; Y(u_0, v_0)=Y_0, \; l(u_0, v_0)=l_0.$$
Furthermore, $(u,v)$ are principal parameters for the surface $S$.
\end{thm}

\section{Natural principal parameters on space-like Weingarten surfaces}

The main object of our considerations are Weingarten space-like surfaces.
We use the usual definition for these surfaces in terms of the principal curvatures
$\nu_1$ and $\nu_2$.
\vskip 1mm

\begin{defn}
A space-like surface $S:\,z=z(u,v)\;(u,v) \in \mathcal D$ with principal curvatures
$\nu_1$ and $\nu_2$ is said to be a {\it space-like Weingarten surface} (or shortly
a {\it space-like W-surface}) if there exist two differentiable functions
$f,\;g$ of one variable, defined in the interval $\mathcal{I} \subseteq {\R}$ and a
differentiable function $\nu = \nu(u,v)\in \mathcal{I}$ defined in $\mathcal D$, such that
$f'\,g'\neq 0, \; f-g\neq 0,$
and the principal curvatures of $S$ at every point are given by
$\nu_1 =f(\nu), \; \nu_2=g(\nu).$
\end{defn}

The following functions are essential in the theory of space-like W-surfaces:
$$I:=\int_{\nu_0}^{\nu} \frac{f'(\nu)\,d \nu}{f(\nu)-g(\nu)}\,,\qquad
J:=\int_{\nu_0}^{\nu} \frac{g'(\nu)\,d \nu}{g(\nu)-f(\nu)}\,.\leqno (3.1)$$
\vskip 2mm
The next statement gives the property of  space-like Weingarten
surfaces, which allows us to introduce special principal parameters
on such surfaces (cf \cite{GM1}).

\begin{lem}\label {L:7.1} Let $S:\; z=z(u,v),\; (u,v)\in \mathcal{D}$
be a  space-like W-surface parameterized by principal parameters.
Then the function
$$\lambda = \sqrt E \exp\left(\int_{\nu_0}^{\nu} \frac{f'\,d \nu}{f-g}\right)
=\sqrt E\, e^{I}$$
does not depend on $v$, while the function
$$\mu = \sqrt G \exp \left(\int_{\nu_0}^{\nu} \frac{-g'\,d \nu}{f-g}\right)
=\sqrt G\, e^{J}$$
does not depend on $u$.
\end{lem}
\vskip 2mm

We define natural principal parameters on a  space-like W-surface as follows:

\begin{defn}
Let $S:\;z=z(u,v),\; (u,v)\in \mathcal{D}$ be a  space-like
W-surface parameterized by principal parameters. The
parameters $(u, v)$ are said to be \emph{natural} principal, if
the functions $\lambda(u)$ and $\mu(v)$ from Lemma \ref{L:7.1} are
constants.
\end{defn}

\begin{prop}\label {P:5.4}
Any  space-like Weingarten surface admits locally natural principal parameters.
\end{prop}

\emph{Proof}: Let $S:\; z=z(u,v),\; (u,v)\in \mathcal{D}$ be a
 space-like W-surface,
parameterized by principal parameters. Then $\nu_1=f(\nu), \; \nu_2=g(\nu), \;
\nu=\nu(u,v)$ for some differentiable functions $f(\nu), \; g(\nu)$ and $\nu(u,v)$
satisfying the condition
$$(f(\nu)-g(\nu))\, f'(\nu)\,g'(\nu) \neq 0, \quad (u,v) \in \mathcal D.$$

Let $\mathfrak{a}={\rm const} \neq 0, \; \mathfrak{b}={\rm const} \neq 0$, \; $(u_0, v_0)
\in \mathcal D$ and $\nu_0=\nu(u_0,v_0)$. We change the parameters
$(u, v)\in \mathcal{D}$ with $(\bar u, \bar v)\in \bar{\mathcal{D}}$
by the formulas
$$\left|\begin{array}{l}
\displaystyle{\bar u=\mathfrak{a}\int_{u_0}^u \sqrt E\,e^I\,du
 +\overline{u}_0 , \quad \bar u_0 = {\rm const},}\\
[4mm]
\displaystyle{\bar v=\mathfrak{b}\int_{v_0}^v\sqrt G \,e^J\,dv
+\overline{v}_0, \quad \bar v_0={\rm const}}.
\end{array}\right.$$
According to Lemma \ref{L:7.1} it follows that $(\bar u, \bar v)$
are again principal parameters and
$$\bar E=\mathfrak{a}^{-2}\,e^{-2I},\quad
\bar G=\mathfrak{b}^{-2}\,e^{-2J}.\leqno (3.2)$$
Then for the functions from Lemma
\ref{L:7.1} we find
$ \lambda(\bar u) =|\mathfrak{a}|^{-1}, \quad  \mu(\bar v) =|\mathfrak{b}|^{-1}.$

Furthermore $\mathfrak{a}^{-2}=\bar E(u_0,v_0), \; \mathfrak{b}^{-2} = \bar G(u_0,v_0)$.
{\qed} \vskip 2mm

\begin{cor}
The principal parameters on a given space-like Weingarten surface are
natural principal if and only if
$$\sqrt{EG}(\nu_1-\nu_2)={\rm const}\neq 0.\leqno (3.3)$$
\end{cor}

Further we assume that the space-like W-surface
$S: z=z(u,v)$, $(u,v) \in \mathcal D$ under consideration, is parameterized
by natural principal parameters $(u, v)$. It follows from the above
proposition that the coefficients $E$ and $G$  (consequently $L$ and
$N$) are expressed by the invariants of the surface.

As an immediate consequence from Proposition \ref{P:5.4} we get

\begin{cor}
Let $S$ be a  space-like W-surface parameterized by
natural principal parameters $(u, v)$. Then any natural
principal parameters $(\tilde u, \tilde v)$ on $S$ are determined by
$(u, v)$ up to an affine transformation of one of the types
$$\begin{array}{ll}
\left|\begin{array}{lcc}
\tilde u =\, a_{11} \, u  & & +\; b_1,\\
[2mm] \tilde v = & a_{22} \, v &+\; b_2, \end{array} \quad
a_{11}a_{22} \neq 0,\right.&
 \quad\left|\begin{array}{lcc}
\tilde u = &a_{12} \, v & + \;c_1,\\
[2mm] \tilde v =\, a_{21} \, u & &+\;c_2, \end{array} \quad
a_{12}a_{21}\neq 0,\right.
\end{array}$$
where $a_{ij}, \, b_i, \, c_i; \, i,j=1,2 $
are constants.
\end{cor}

\subsection{Strongly regular space-like W-surfaces.}
First we consider strongly regular space-like W-surfaces, i.e. W-surfaces, satisfying the
condition
$$\nu_u(u,v)\nu_v(u,v)\neq 0, \qquad (u,v) \in \mathcal D.$$

Our main theorem for such surfaces states

\begin{thm}\label {T:7.3}
Given two differentiable functions $f(\nu), \; g(\nu); \; \nu \in
\mathcal{I},$ $f(\nu)-g(\nu)\neq 0,$ $f'(\nu)g'(\nu)\neq 0$ and a
differentiable function $\nu(u,v), \; (u,v) \in {\mathcal D}$
satisfying the conditions
$$\nu_u\,\nu_v\neq 0, \quad \nu(u,v)\in \mathcal{I}.$$

Let $(u_0, v_0) \in \mathcal D, \; \nu_0=\nu(u_0, v_0)$ and $\mathfrak{a}>0, \, \mathfrak{b}> 0$
be two constants. If
$$\mathfrak{a}^2\,e^{2 I}\,(J_{uu}+I_u\,J_u-J_u^2)+
\mathfrak{b}^2\,e^{2 J}(I_{vv}+I_v\,J_v-I_v^2)
=-f(\nu)\,g(\nu),\leqno(3.4)$$
then there exists a unique (up to a motion) strongly regular
space-like W-surface $S:\; z=z(u,v),$ \, $(u,v)\in \mathcal D_0 \subset \mathcal D$
with invariants
$$\begin{array}{ll}
\nu_1=f(\nu), & \nu_2=g(\nu),\\
[2mm]
\gamma_1=\mathfrak{b}\,e^{J}\,(I)_v,& \gamma_2=-\mathfrak{a}\,e^{I}\,(J)_u.
\end{array}\leqno(3.5) $$
Furthermore, $(u,v)$ are natural principal parameters on $S$.
\end{thm}

We also have
$$\mathfrak{a}\,\sqrt{E}= e^{- I}, \quad   \mathfrak{b}\,\sqrt{G}= e^{- J}.$$

Hence, with respect to natural principal parameters each strongly regular space-like
Weingarten surface possesses a {\it natural PDE} (3.4).

\vskip 2mm
\subsection{Space-like W-surfaces with $\gamma_1=0$.}

In this subsection we consider space-like W-surfaces with $\gamma_1=0$
and prove the fundamental theorem of Bonnet type for this class.

Let $S:\; z=z(u,v),\; (u,v)\in \mathcal{D}$ be a  space-like W-surface,
parameterized by natural principal parameters. Then we can assume
$$\sqrt E=\frac{1}{\mathfrak{a}}\,e^I, \qquad \sqrt G=\frac{1}{\mathfrak{b}}\,e^J,$$
where $I$ and $J$ are the functions (3.1) and $\mathfrak{a},\,\mathfrak{b}$
are some positive constants.
We note that under the condition $\gamma_1=0$ it follows that the function
$\nu=\nu(u)$ does not depend on $v$.

Considering the system (2.3), we obtain that the compatibility conditions for
this system reduce to only one - the Gauss equation, which has the form:
$$X(\gamma_2)+\gamma_2^2=f(\nu)\,g(\nu).$$

Thus we obtain the following Bonnet type theorem for space-like W-surfaces
satisfying the condition $\gamma_1=0$:
\begin{thm}
Given two differentiable functions $f(\nu), \; g(\nu); \; \nu \in
\mathcal{I},$ $f(\nu)-g(\nu)\neq 0,$ $f'(\nu)g'(\nu)\neq 0$ and a
differentiable function $\nu(u,v)=\nu(u), \; (u,v) \in {\mathcal D}$
satisfying the conditions
$$\nu_u\neq 0, \quad \nu(u,v)\in \mathcal{I}.$$

Let $(u_0, v_0) \in \mathcal D, \; \nu_0=\nu(u_0, v_0)$ and $\mathfrak{a}>0$
be a constant. If
$$\mathfrak{a}^2\,e^{2 I}\,(J_{uu}+I_u\,J_u-J_u^2)=-f(\nu)\,g(\nu),\leqno(3.6)$$
then there exists a unique (up to a motion) space-like W-surface
$S:\; z=z(u,v),$ \, $(u,v)\in \mathcal D_0 \subset \mathcal D$
with invariants
$$\begin{array}{ll}
\nu_1=f(\nu), & \nu_2=g(\nu), \\
[2mm]
\gamma_1=0, & \gamma_2=-\mathfrak{a}\,e^{I}\,(J)_u.
\end{array}\leqno(3.7) $$
Furthermore, $(u,v)$ are natural principal parameters on $S$.
\end{thm}

Hence, with respect to natural principal parameters each space-like
Weingarten surface with $\gamma_1=0$ possesses a {\it natural ODE} (3.6).
\vskip 2mm

\section{Parallel space-like surfaces in Minkowski space and their natural PDE's}

Let $S: \; z=z(u,v),\; (u,v)\in \mathcal{D}$ be a space-like surface, parameterized by
principal parameters and $l(u,v),\;l^2=-1$ be the unit normal vector field of $S$.
The parallel surfaces of $S$ are given by
$$\bar S(a): \; \bar z(u,v)= z(u,v) + a\,l(u,v), \quad a={\rm const} \neq 0,
\quad (u,v)\in \mathcal{D}.\leqno(4.1)$$

We call the family $\{\bar S(a), \; a={\rm const} \neq 0\}$ the \emph{parallel
family} of $S$.

Taking into account (4.1), we find
$$\bar z_{u}=(1-a\,\nu_1)\, z_u,\quad\bar z_{v}=(1-a\,\nu_2)\, z_v.\leqno(4.2)$$

Excluding the points, where $(1-a\,\nu_1)(1-a\,\nu_2)=0$, we obtain that the corresponding
unit normal vector fields $\bar l$ to $\bar S(a)$ and $l$ to $S$ satisfy the equality
$\bar l=\varepsilon\,l,$ where $\varepsilon:= {\rm sign}\,(1-a\,\nu_1)(1-a\,\nu_2).$
Hence, the parallel surfaces $\bar S(a)$ of a space-like surface $S$ are also space-like
surfaces.

The relations between the principal curvatures $\nu_1(u,v)$,
$\nu_2(u,v)$ of $S$ and $\bar \nu_1(u,v)$,  $\bar \nu_2(u,v)$ of its
parallel space-like surface $\bar S(a)$ are
$$\bar \nu_1=\varepsilon\,\frac{\nu_1}{1-a\,\nu_1}\,, \quad \bar \nu_2=\varepsilon\,
\frac{\nu_2}{1-a\,\nu_2};\quad \nu_1=\frac{\varepsilon\,\bar\nu_1}{1+a\,\varepsilon\,\bar\nu_1}\,,
\quad \nu_2=\frac{\varepsilon\,\bar\nu_2}{1+a\,\varepsilon\,\bar\nu_2}.\leqno (4.3)$$

Let $K'=\nu_1\,\nu_2,\; \displaystyle{H=\frac{1}{2}(\nu_2+\nu_2),\;H'=\frac{1}{2}(\nu_2-\nu_2)}$
be the three invariants of the space-like surface $S$. The equalities (4.3) imply the  relations
between the invariants $\bar K'$, $\bar H$ and $\bar H'$ of $\bar S(a)$ and the corresponding invariants of $S$:
$$K'=\frac{\bar K'}{1+2 a\,\varepsilon \bar H + a^2 \bar K'}\,, \quad
H=\frac{\varepsilon\,\bar H+a \bar K'}{1+2 a\,\varepsilon \bar H + a^2 \bar K'}\,,\quad
H'=\frac{\varepsilon\,\bar H'}{1+2 a\,\varepsilon\, \bar H + a^2 \bar K'}\,.\leqno (4.4)$$
\vskip 2mm

Now let $S: \; z=z(u,v),\; (u,v)\in \mathcal{D}$ be a space-like Weingarten surface with
Weingarten functions $f(\nu)$ and $g(\nu)$. We suppose that $(u,v)$ are
natural principal parameters for $S$.
We show that $(u,v)$ are also natural principal parameters
for any parallel space-like surface $\bar S(a)$.
\begin{prop}
The natural principal parameters $(u,v)$ of a given space-like W-surface $S$ are natural
principal parameters for all parallel space-like surfaces $\bar S(a),\; a={\rm const}
\neq 0$ of $S$.
\end{prop}
{\it Proof:} Let $(u,v)\in \mathcal{D}$ be natural principal parameters for $S$,
$(u_0,v_0)$ be a fixed point in $\mathcal D$ and $\nu_0=\nu(u_0,v_0)$.
The coefficients $E$ and $G$ of the first fundamental form of $S$ are given by (3.2).
The corresponding coefficients $\bar E$ and $\bar G$ of $\bar S(a)$ in view of
(4.2) are
$$\bar E=(1-a\,\nu_1)^2\,E,\quad \bar G=(1-a\,\nu_2)^2\,G.\leqno(4.5)$$
Equalities (4.3) imply that $\bar S(a)$ is again a Weingarten surface with Weingarten
functions
$$\bar\nu_1(u,v)=\bar f(\nu)=\frac{\varepsilon f(\nu)}{1-af(\nu)}\,,
\quad\bar\nu_2(u,v)=\bar g(\nu)=\frac {\varepsilon g(\nu)}{1-ag(\nu)}. \leqno(4.6)$$
Using (4.6), we compute
$$\bar f-\bar g=\frac{\varepsilon(f-g)}{(1-a\,f)(1-a\,g)}\,,$$
which shows that ${\rm sign} \,(\bar f-\bar g)={\rm sign}\,(f-g)$.

Further, we denote by $f_0:=f(\nu_0), \; g_0:=g(\nu_0)$ and taking into account
(3.3) and (4.5), we compute
$$\sqrt{\bar E\,\bar G}\,(\bar f-\bar g)=\sqrt{E\, G}\,( f- g)=
{\rm const},$$
which proves the assertion. \qed
\vskip 2mm
Using the above statement, we prove the following theorem.
\begin{thm}
The natural PDE of a given space-like W-surface $S$ is the natural PDE of any parallel
space-like surface $\bar S(a),\; a={\rm const}\neq 0$, of $S$.
\end{thm}
{\it Proof}. We have to express the equation (3.4) (resp. (3.6)) in terms of the
Weingarten functions of the parallel space-like surface $\bar S(a)$.
Using (3.1) and (4.6), we compute
$$\bar I= I-\ln\frac{1-a\,f}{1-a\,f_0},\quad
\bar J= J-\ln\frac{1-a\,g}{1-a\,g_0}.$$

Putting
$$\bar E_0=(1-a\,\nu_{1}(u_0,v_0))^2\,E_0=\mathfrak{a}^{-2}\,(1-a\,f_0)^2
=:\bar {\mathfrak{a}}^{-2},$$
$$\bar G_0=(1-a\,\nu_{2}(u_0,v_0))^2\,G_0=\mathfrak{b}^{-2}\,(1-a\,g_0)^2
=:\bar {\mathfrak{b}}^{-2},$$
we obtain
$$\begin{array}{l}
\bar{\mathfrak{a}}^2\,e^{2 \bar I}\,(\bar{J}_{uu}+\bar{I}_u\,\bar{J}_u-\bar{J}_u^2)+
\bar{\mathfrak{b}}^2\,e^{2 \bar{J}}(\bar{I}_{vv}+\bar{I}_v\,\bar{J}_v-\bar{I}_v^2)
+\bar{f}(\nu)\,\bar{g}(\nu)\\
[2mm]
=\mathfrak{a}^2\,e^{2 I}\,(J_{uu}+I_u\,J_u-J_u^2)+
\mathfrak{b}^2\,e^{2 J}(I_{vv}+I_v\,J_v-I_v^2)
+f(\nu)\,g(\nu).
\end{array}$$

Hence, the natural PDE of $\bar S(a)$ in terms of the Weingarten functions
$\bar f(\nu)$, $\bar g(\nu)$ coincides with the natural PDE of $S$ in terms of the
Weingarten functions $f(\nu)$ and $g(\nu)$. \qed

\section{Linear fractional space-like Weingarten surfaces}

In this section we study special classes of space-like Weingarten surfaces
with respect to the Weingarten functions $f$ and $g$.
\vskip 1mm

A space-like W-surface with principal curvatures $\nu_1$ and
$\nu_2$ is said to be {\it linear fractional} if
$$\nu_1=\frac{A\nu_2+B}{C\nu_2+D}\,, \quad A, B, C, D - {\rm constants},\quad BC-AD\neq 0.\leqno(5.1)$$

We exclude the case $A=D,\,B=C=0$, which characterizes the umbilical points.

\begin{lem}
Any linear fractional space-like Weingarten surface determined by $(5.1)$ is
a surface whose invariants $K', H, H'$ satisfy the linear relation
$$\delta K'=\alpha \,H+\beta\,H'+\gamma, \quad \alpha, \beta, \gamma, \delta - {\rm constants},
\quad \alpha^2-\beta^2+4\gamma\delta\neq 0\leqno(5.2)$$
and vice versa.
\end{lem}

The relations between the constants $A, B, C, D$ in (5.1) and $\alpha, \beta, \gamma, \delta$
in (5.2) are given by the equalities:
$$\alpha=A-D, \quad \beta = -(A+D), \quad \gamma=B,
\quad \delta=C.\leqno(5.3)$$

We denote by $\mathcal K$ the class of all space-like surfaces, free
of umbilical points, whose curvatures satisfy (5.2) or equivalently
(5.1).

The aim of our study is to classify all natural PDE's of the surfaces from the class
$\mathcal K$.

The parallelism between two surfaces
given by (4.1) is an equivalence relation. On the other hand, Theorem 4.2 shows that
the surfaces from an equivalence class have one and the same natural PDE. Hence, it is
sufficient to find the natural PDE's of the equivalence classes. For any equivalence class,
we use a special representative, which we call \emph{a basic class}. Thus the classification
of the natural PDE's of the surfaces in the class $\mathcal K$ reduces to the
natural PDE's of the basic classes.

In view of Theorem 4.2, we prove the following classification
Theorem.

\begin{thm}
 Up to similarity, the space-like surfaces in Minkowski space, whose
curvatures $K'$, $H$ and $H'$ satisfy the linear relation
$$\delta K' = \alpha H + \beta H' + \gamma, \quad \alpha, \beta, \gamma, \delta -
{\rm constants}; \quad \alpha^2-\beta^2 + 4 \gamma\delta \neq 0,$$
are described by the natural PDE's of the following basic surfaces:
\vskip 3mm
\begin{enumerate}
\item
$H=0:\quad  \nu=-e^{\lambda},\quad \Delta \lambda =  e^{\lambda},$\\
\vskip 3mm
\item
$H=\displaystyle{\frac{1}{2}:\quad \nu=\frac{1}{2}(1-e^{\lambda}),\quad
\Delta \lambda = \sinh{\lambda}},$\\
\vskip 3mm
\item
$H'=1: \quad \Delta^*(e^{\nu})=2\,\nu\,(\nu+2),$\\
\vskip 3mm
\item
$\left|\begin{array}{l}
H=p\,H'\\
[1mm]
p^2>1\end{array}\right.:\qquad \Delta^*(\nu^{p})=
\displaystyle{2\,\frac{p\,(p+1)}{(p-1)^2}\,\nu},$\\
\vskip 3mm
\item
$\left|\begin{array}{l}
H=p\,H'\\
[1mm]
p^2<1\end{array}\right.:\qquad \bar\Delta^*(\nu^{p})=
\displaystyle{2\,\frac{p\,(p+1)}{(p-1)^2}\,\nu},$\\
\vskip 3mm
\item
$\left|\begin{array}{l}
H=p\,H'+1 \\
[1mm]
p^2>1\end{array}\right.:
\quad  \displaystyle{\nu=\frac{(p-1)\,\lambda+2}{2},\quad
\Delta^*(\lambda^{p})= \frac{p\,((p-1)\lambda+2)((p+1)\lambda+2)}
{2\,(p-1)\,\lambda}},$\\
\vskip 3mm
\item
$\left|\begin{array}{l}
H=p\,H'+1 \\
[1mm]
p^2<1\end{array}\right.:
\quad  \displaystyle{\nu=\frac{(p-1)\,\lambda+2}{2},\quad
\bar\Delta^*(\lambda^{p})= \frac{p\,((p-1)\lambda+2)((p+1)\lambda+2)}
{2\,(p-1)\,\lambda}},$\\
\vskip 3mm
\item
$K'=-1:\quad \nu=\tan\lambda,\quad \bar\Delta \lambda = -\sin \lambda,$\\
\vskip 3mm
\item
$K'=2\,H': \quad \displaystyle{\nu=\frac{\lambda-4}{\lambda-2},\quad
\Delta^*(e^{\lambda})=2},$\\
\vskip 3mm
\item
$\left|\begin{array}{l}
K'=p\,H'-q\\
[1mm]
p\neq 0,\,q>0\end{array}\right.:
\quad \left |
\begin{array}{l}
\displaystyle{\nu=\lambda+\frac{p}{2}},\\
[3mm]
\displaystyle{\mathcal{I}=\frac{1}{\sqrt{q}}\,\arctan\frac{\lambda}{\sqrt{q}}}
\end{array} \right.,\quad
 \displaystyle{\Delta^*(e^{p\,\mathcal{I}}) =\frac{p\,q}{2}\,
\frac{\lambda\,\left(p\,\lambda-2\,q\right)}{\lambda^2+q}}.$
\end{enumerate}
\end{thm}
\vskip 2mm

{\it Proof:}
According to the constant $C$ in (5.1), the linear fractional space-like W-surfaces
are divided into two classes:
linear fractional space-like W-surfaces, determined by the condition $C=0$
and linear fractional space-like W-surfaces, determined by
the condition $C\neq 0$.

{\bf I.}\; Linear fractional space-like Weingarten surfaces with
$C=0$. This class is determined by the equality
$$\alpha\,H+\beta\,H'+\gamma=0,\quad (\alpha,\gamma)\neq (0,0),\quad
\alpha^2-\beta^2\neq 0.\leqno(5.4)$$

For the invariants of the space-like parallel surface $\bar S(a)$ of $S$, because of (4.4),
we get the relation
$$\varepsilon\,(\alpha+2\,a\,\gamma)\,\bar H+\varepsilon\,\beta\,\bar H'+\gamma
=-a\,(\alpha+a\,\gamma)\,\bar K'.\leqno (5.5)$$

Let $\;\eta:= {\rm sign}\,(\alpha^2-\beta^2).$ Each time choosing
appropriate values for the constants $\mathfrak{a}$, $\mathfrak{b}$
and $\nu_0$ in (3.4), we consider the following subclasses and their
natural PDE's:
\vskip 2mm

\begin{itemize}
\item[1)] $\alpha=0,\;\beta\neq 0,\;\gamma\neq 0$. Assuming that $\gamma=1$, the relation
(5.4) becomes
$$\beta\,H'+1=0. $$

The natural PDE for these W-surfaces is
$$(e^{\beta\,\nu})_{vv}
+(e^{-\beta\,\nu})_{uu}=
\frac{2}{\beta}\,\nu\,(\beta\,\nu-2).\leqno (5.6)$$

Up to similarities these W-surfaces are generated by the basic class $H'=1$
with the natural PDE
$$(e^{\nu})_{uu}
+(e^{-\nu})_{vv}=
2\,\nu\,(\nu+2),\leqno (5.6^*)$$
which is the case $(3)$ in the statement of the theorem.
\vskip 2mm

\item[2)]
 $\displaystyle{\alpha\neq 0,\;\gamma=0}$. Assuming that $\alpha=1$, the relation
 (5.4) becomes
$$H+\beta\,H'=0.$$
\begin{itemize}
\item[2.1)] $\beta\neq 0,\;\eta=-1\; (\beta^2-1>0).$
Choosing
$\displaystyle{\mathfrak{b}^2\,\frac{\beta-1}{\beta+1}\,\nu_0^{-(\beta+1)}=1,\;
\mathfrak{a}^2\,\nu_0^{\beta-1}=1}$, the natural PDE becomes
$$\left(\nu^{\beta}\right)_{vv}+
\left(\nu^{-\beta}\right)_{uu}=
2\,\frac{\beta(\beta-1)}{(\beta+1)^2}\,\nu,\leqno(5.7)$$
which is the case $(4)$ in the statement of the theorem.
\vskip 2mm

\item[2.2)] $\beta \neq 0,\;\eta=1\; (\beta^2-1<0)$. Choosing
$\displaystyle{\mathfrak{b}^2\,\frac{\beta-1}{\beta+1}\,\nu_0^{-(\beta+1)}=-1,\;
\mathfrak{a}^2\,\nu_0^{\beta-1}=1}$, the natural PDE becomes
$$\left(\nu^{\beta}\right)_{vv}-
\left(\nu^{-\beta}\right)_{uu}=
-2\,\frac{\beta(\beta-1)}{(\beta+1)^2}\,\nu\leqno (5.8)$$
which is the case $(5)$ in the statement of the theorem.
\vskip 2mm

\item[2.3)] $\beta = 0$. Putting  $\;\nu =e^{\lambda}$, we get the natural
PDE for space-like surfaces with $H=0$:
$$\Delta \lambda = e^{\lambda},\leqno (5.9)$$
which is the case $(1)$ in the statement of the theorem.
\end{itemize}
\vskip 2mm

\item[3)] $\alpha\neq 0,\;\beta=0,\; \gamma\neq 0.$ Assuming that $\alpha=1$,
the relation (5.4) becomes
$$H+\gamma=0.$$
Putting  $\;\displaystyle{| H|\,e^{\lambda}:=H-\nu}=H'>0$,
we get the one-parameter system of natural PDE's for CMC space-like surfaces
with $H=-\gamma$:
$$\Delta \lambda= 2\,|H|\,\sinh \lambda.\leqno(5.10)$$

Up to similarities these W-surfaces are generated by the basic class $|H|=\frac{1}{2}$
with the natural PDE
$$\Delta \lambda= \sinh \lambda,\leqno(5.10^*)$$
which is the case $(2)$ in the statement of the theorem.
\vskip 2mm

\item[4)] $ \alpha\neq 0,\;\beta\neq 0,\;\gamma\neq 0.$

Assuming that $\alpha=1$ we have
$$H+\beta\,H'+\gamma=0,\quad \beta^2-1\neq 0.$$
Let $\displaystyle{\lambda:=2\,H'=\frac{-2}{\beta+1}\,(\nu+\gamma)}>0$.
\vskip 2mm

\begin{itemize}
\item[4.1)] If $\eta=-1 \; (\beta^2-1>0)$,
and choosing
$$\mathfrak{b}^2=\frac{\beta+1}{\beta-1}\,\left(\frac{-2}{\beta+1}
(\nu_0+\gamma)\right)^{\beta+1},\;
\mathfrak{a}^2\,=\left(\frac{-2}{\beta+1}(\nu_0+\gamma)\right)^{-(\beta-1)},$$

the natural PDE becomes
$$\left(\lambda^{\beta}\right)_{vv}
+\left(\lambda^{-\beta}\right)_{uu}=
\frac{\beta}{2\,(\beta+1)}\,
\frac{((\beta+1)\lambda+2\,\gamma)((\beta-1)\lambda+2\,\gamma)}{\lambda}\,.
\leqno (5.11)$$

Up to similarities these W-surfaces are generated by the basic class\\
$H=\beta\,H'+1,\;\beta^2>1$ with the natural PDE
$$\left(\lambda^{\beta}\right)_{uu}
+\left(\lambda^{-\beta}\right)_{vv}=
\frac{\beta}{2\,(\beta-1)}\,
\frac{((\beta+1)\lambda+2)((\beta-1)\lambda+2)}{\lambda}\,,
\leqno (5.11^*)$$
which is the case $(6)$ in the statement of the theorem.
\vskip 2mm

\item[4.2)] If $\eta=1\; (\beta^2-1<0)$, and choosing
$$\mathfrak{b}^2=-\frac{\beta+1}{\beta-1}\,\left(\frac{-2}{\beta+1}
(\nu_0 + \gamma)\right)^{\beta+1},\;
\mathfrak{a}^2\,=\left(\frac{-2}{\beta+1}(\nu_0 + \gamma)\right)^{-(\beta-1)},$$

the natural PDE becomes
$$\left(\lambda^{\beta}\right)_{vv}
-\left(\lambda^{-\beta}\right)_{uu}=
\frac{-\beta}{2\,(\beta+1)}\,
\frac{((\beta+1)\lambda+2\,\gamma)((\beta-1)\lambda+2\,\gamma)}{\lambda}\,.
\leqno (5.12)$$

Up to similarities these W-surfaces are generated by the basic class\\
$H=\beta\,H'+1,\;\beta^2<1$ with the natural PDE
$$\left(\lambda^{\beta}\right)_{uu}
-\left(\lambda^{-\beta}\right)_{vv}=
\frac{\beta}{2\,(\beta-1)}\,
\frac{((\beta+1)\lambda+2)((\beta-1)\lambda+2)}{\lambda}\,,
\leqno (5.12^*)$$
which is the case $(7)$ in the statement of the theorem.
\end{itemize}
\end{itemize}

{\bf II.}\;
Linear fractional space-like Weingarten surfaces with $C\neq 0$.
Let $C=1$. The equality (5.2) gets the form
$$K'=\alpha\,H+\beta\,H'+\gamma, \quad \alpha^2-\beta^2+4\,\gamma\neq 0.\leqno(5.13)$$

The corresponding relation for the parallel surface $\bar S(a)$ is
$$\varepsilon(\alpha+2\,a\,\gamma)\,\bar H+\varepsilon\,\beta\,\bar H'+\gamma=
(1-a\,\alpha-a^2\,\gamma)\,\bar K'.\leqno(5.14)$$

Each time choosing appropriate values for the constants
$\mathfrak{a}$, $\mathfrak{b}$ and $\nu_0$ in (3.4), we consider the following subclasses
and their natural PDE's:

\vskip 3mm
\begin{itemize}
\item[5)] $\alpha = \gamma= 0, \; \beta \neq 0$. The relation (5.13) becomes
$$K'=\beta H' \quad \iff \quad \rho_1-\rho_2=-\frac{2}{\beta}\,,$$
where $\displaystyle{\rho_1=\frac{1}{\nu_1},\;\rho_2=\frac{1}{\nu_2}}$ are the principal
radii of curvature of $S$.

Putting $\lambda:=\displaystyle{4\,\frac{\nu-\beta}{2\,\nu-\beta}}$,
and choosing $\nu_0=\beta$,
the natural PDE of these space-like surfaces gets the form
$$\left(e^\lambda\right)_{uu} +\left(e^{-\lambda}\right)_{vv}-\frac{\beta^4}{8}=0.
\leqno (5.15)$$

Up to similarities these W-surfaces are generated by the basic class
$K'=2\,H'$ with the natural PDE
$$\left(e^\lambda\right)_{uu} +\left(e^{-\lambda}\right)_{vv}-2=0,
\leqno (5.15^*)$$
which is the case $(9)$ in the statement of the theorem.

\vskip 2mm
\begin{rem}
A. Ribaucour \cite{Rib} has proved that
{\it a necessary condition for the curvature lines
of the first and second focal surfaces of $S$ to  correspond to each other resp. to
a conjugate parametric lines on $S$ is $\rho_1-\rho_2={\rm const}$ resp.
$\rho_1\,\rho_2={\rm const}$}.

Von Lilienthal \cite{vL1} (cf \cite{vL2, Bi1, Bi2, E})
has proved in $\mathbb{R}^3$
that a surface with
$\rho_1-\rho_2=1/R,\; R={\rm const}\neq 0,$ has  first and second focal surfaces
 of constant Gauss curvature $-R^2$ and vice versa.

In $\mathbb{R}^3_1$ one can prove in a similar way the corresponding property:
The first and second focal surfaces
of a space-like surface with $K'=\beta\,H'$ are
time-like  of constant Gauss curvature $\beta^2/4$, and vice versa.
\end{rem}
\vskip 2mm

\item [6)] $(\alpha, \gamma)\neq (0,0), \; \alpha^2+4\gamma\geq 0$. The relation (5.14)
implies that there exists a space-like surface $\bar S(a)$, parallel to $S$, which satisfies the
relation (5.4). Hence the natural PDE of $S$ is one of the PDE's (5.6) - (5.12).
\vskip 3mm

\item[7)] $\alpha^2+4\,\gamma <0$. It follows that $\gamma <0$. The relation (5.14)
implies that there exists a surface $\bar S(a)$ parallel to $S$, which satisfies the relation
$$K'=\beta H'+\gamma.\leqno (5.16)$$
\begin{itemize}
\item[7.1)] $\beta = 0$. The relation (5.16) becomes
$K'=\gamma<0,$ i.e. $\bar S$ is of constant negative extrinsic sectional
curvature $\gamma$ (of constant positive intrinsic sectional
curvature $K=-\gamma$). Putting
$\displaystyle{\lambda:=2\,\arctan\frac{\nu}{\sqrt{-\gamma}}}$, we
get the natural PDE of this surface
$$\bar\Delta \lambda=-K^2\,\sin \lambda.
\leqno (5.17)$$

Up to similarities these W-surfaces are generated by the basic class
$K'=-1$ ($K=1$) with the natural PDE
$$\bar\Delta \lambda=-\sin \lambda,\leqno (5.17^*)$$
which is the case $(8)$ in the statement of the theorem.
\vskip 2mm

\item[7.2)] $\beta \neq 0$, $\gamma< 0$. Choosing
$\displaystyle{\nu_0=\frac{\beta}{2}}$,
 the natural PDE of $S$ becomes
$$(\exp{(\beta\,\mathcal{I})})_{uu}+
(\exp{(-\beta\,\mathcal{I})})_{vv}
=-\frac{\beta\,\gamma}{2}\,
\frac{\lambda\,\left(\beta\,\lambda+2\,\gamma\right)}{\lambda^2-\gamma}\,,\leqno(5.18)$$
where
$$\mathcal{I}=\frac{1}{\sqrt{-\gamma}}
\,\arctan\frac{\lambda}{\sqrt{-\gamma}},\quad
\lambda:=\nu-\frac{\beta}{2},$$
which is the case $(10)$ in the statement of the theorem.
\end{itemize}

\end{itemize}

$\hfill\square$

\end{document}